\documentclass[letterpaper, 10 pt, journal, twoside]{IEEEtran}
\usepackage{amssymb,amsmath}
\usepackage{amsfonts}
\usepackage[english]{babel}
\usepackage[latin1]{inputenc}
\usepackage{fancyhdr}
\usepackage{yfonts}
\usepackage{mathrsfs}
\usepackage[dvips]{graphicx}
\usepackage{enumerate}
\usepackage{xspace}
\usepackage{cite}

\usepackage{psfrag}
\usepackage{url}

\usepackage[dvipsnames]{xcolor}
\usepackage{xcolor}
\usepackage{lipsum}

\usepackage{verbatim}

\providecommand{\com[2]}{\begin{tt}[#1: #2]\end{tt}}
\usepackage{color}

\newtheorem{theorem}{Theorem}
\newtheorem{lemma}[theorem]{Lemma}

\providecommand{\prt}[1]{\left( #1 \right)}

\providecommand{\abs[1]}{\left|#1\right|}
\providecommand{\norm[1]}{\abs{\abs{#1}}}%

\providecommand{\1}{\textbf{1}}
\providecommand{\0}{\textbf{0}}\providecommand{\quotes[1]}{``#1"}

\def\0{{\bf 0}}
\def\R{\mathbb{R}}

\providecommand{\RR}{\mathbb{R}}

\begin{document}

\title{Trajectory convergence from coordinate-wise decrease of quadratic energy functions,\\  and  applications to platoons}

\author{Julien M. Hendrickx, Bal\'azs Gerencs\'er and Baris Fidan
\thanks{J. Hendrickx is with the ICTEAM institute, UCLouvain, Belgium, and the CISE, Boston University, USA. {\tt julien.hendrickx @uclouvain.be} His work is supported by ``Communaut\'e fran\c{c}aise de Belgique -
Actions de Recherche Concert\'ees'' and a WBI.World excellence fellowship.}
\thanks{B. Gerencs\'er is with MTA Alfr\'ed R\'enyi Institute of Mathematics, Budapest, Hungary and E\"otv\"os Lor\'and University, Department of Probability and Statistics, Budapest, Hungary.
{\tt gerencser.balazs@renyi.mta.hu}
He was supported by NKFIH (National Research, Development and Innovation Office) Grants PD 121107 and KH 126505.}
\thanks{B. Fidan is with the Mechanical and Mechatronics Engineering Department, University of Waterloo, ON, Canada. {\tt fidan@uwaterloo.ca} His work is supported by Natural Sciences and Engineering Research Council (NSERC) of Canada under Discovery Grant 116806.}
}

\maketitle

\thispagestyle{empty}

\begin{abstract}
We consider trajectories where the sign of the derivative of each entry is opposite to that of the corresponding entry in the gradient of an energy function. We show that this condition guarantees convergence when the energy function is quadratic and positive definite and partly extend that result to some classes of positive semi-definite quadratic functions including those defined using a graph Laplacian.
We show how this condition allows establishing the convergence of a platoon application in which it naturally appears, due to deadzones in the control laws designed to avoid instabilities caused by inconsistent measurements of the same distance by different agents.
\end{abstract}

\begin{IEEEkeywords}
Lyapunov methods; Agents-based systems; Autonomous vehicles
\end{IEEEkeywords}

\section{Introduction}

We consider trajectories $x(t):\RR_+ \to \RR^n$ that satisfy the following coordinate-wise condition
\begin{equation}\label{eq:general_condition}
\dot x_i \cdot (\nabla V(x))_i \leq 0,  \hspace{.3cm} \forall i=1,\dots,n
\end{equation}
for some quadratic energy function $V$, where the $i$th entry of any vector $v$ is denoted by $v_i$. Intuitively, this just requires that when a coordinate of $x$ changes,  this change happens in a way that does not increase the energy function, but there is no requirement about the magnitude of the decrease, if any.  We show that Condition \eqref{eq:general_condition} guarantees convergence of $x$ to a constant vector when $V$ is positive definite, and partly extend this result to classes of positive semi-definite functions. We prove in particular that convergence is still guaranteed when the matrix defining $V$ is a graph Laplacian, so that $V$ is a measure of the \quotes{disagreement} between the $x_i$.

Condition \eqref{eq:general_condition} appears naturally in certain multi-agent dynamics, including the platooning problem we will analyze in Section \ref{sec:platoon}, which involves control laws with deadzones to remove potential instabilities resulting from incoherent measurements of the same distance by different agents. And indeed, we came across Condition \eqref{eq:general_condition} when trying to establish the convergence of the system of Section \ref{sec:platoon}.

Classical approaches for establishing convergence based on energy functions rely on variation of the Lyapunov - Kraskowski - LaSalle theorems \cite{khalil2002nonlinear,ignatyev2013construction}. These approaches apply to unforced
dynamical systems of the form $\dot x = f(x,t)$, where the vector field $f$ often satisfies some (uniform) continuity condition \cite{damak2014global,balan2005extension}.
For example, LaSalle theorem guarantees (under certain conditions) the convergence of solutions of $\dot x = f(x)$  to an invariant set, but not necessarily to a point, provided that  $V(x(t))$ is nonincreasing everywhere \cite{lasalle1960some}.
Convergence to a single point is only guaranteed under additional conditions, such as $\frac{d}{dt}V(x(t))$ being sufficiently negative, by one of the original Lyapunov theorems. Another result for time-varying systems guarantees convergence to 0 if  $\frac{d}{dt}V(x(t))\leq 0 $ and is not identically 0 on any trajectory other than that staying at 0 \cite{sreedhar1970concerning}. For a survey on various cases of unforced systems we direct the reader to \cite{sepulchre1997constructive} as a starting point.

Particularly in cyber-physical systems or systems involving discrete computations or events, there may not be a natural way of defining a global evolution of the form $\dot x = f(x,t)$, for instance when external noise or control is present. Therefore $x$ may not contain all the information required to determine $\dot x$ and would thus not qualify as \quotes{the state of the system} in a classical sense. The speed $\dot x$ may indeed depend on various elements related to the history of $x$, communications with other systems, random or arbitrary events, etc.
Certain works on consensus overcome this difficulty by defining a trajectory-dependent equivalent vector field $\tilde f_x(x,t)$ to hide the complexity of the process, i.e. a vector field for which $\dot x = \tilde f_x(x,t)$ holds for that specific trajectory $x$ only\cite{de2013self,hendrickx2011new,jadbabaie2003coordination}. But it can be challenging to define equivalent fields satisfying the global conditions required to apply the classical convergence results, (continuity, decrease of $V$, suitable invariant sets...), and we indeed did not succeed in applying this approach to general trajectories satisfying \eqref{eq:general_condition}. Moreover, we would argue that this is a cumbersome and unnatural step. Extensions of Lyapunov results to differential inclusions could also not be directly applied to \eqref{eq:general_condition} as they require pre-defining the invariant sets to which $x$ would converge \cite{bacciotti1999stability,corts2008discontinuous}. Hence we think it is in many cases relevant to analyze the convergence of trajectories based purely on their properties, and not on those of a vector field or differential inclusion they follow.

Standard trajectory-focused techniques do not allow establishing convergence solely based on \eqref{eq:general_condition}. Observe it implies 
$$
\frac{d}{dt} V(x(t)) = \sum_{i=1}^n (\nabla V(x))_i \dot x_i \leq 0, $$
so that $V(x(t))$ is non-increasing and hence converging, but it is well known that $\frac{d}{dt} V(x(t))\leq 0$ does not imply the convergence of $x$ in general. There is here no guarantee on the decrease rate of $V$, even relative to the magnitude of $\dot x$, as the gradient $\nabla V(x)$ and $\dot x$ can be orthogonal or arbitrarily close to being orthogonal. Hence the total decrease of $V$ cannot be directly bounded relative to the length of the trajectory, which would have guaranteed a finite length of the trajectory. We can thus a priori not exclude that $x$ would keep varying while approaching a level set $\{x:V(x)=c\}$.
Moreover, these level sets are not necessarily compact when $V$ is only a positive semi-definite quadratic function.

In the following section we establish our new convergence result based on Condition \ref{eq:general_condition}. Afterwards, we demonstrate its application for a platoon formation problem.

\section{Convergence result}\label{sec:convergence}

For the simplicity of exposition, we state our main convergence result for quadratic functions of the form $x^TAx$ and particularize condition \eqref{eq:general_condition} to these functions, but extension to general quadratic functions is immediate by applying a constant offset $x'= x-a$ for some vector $a$.

\begin{theorem}\label{thm:elementwise}
  Let $A\in \RR^{n\times n}$ be a symmetric positive semi-definite matrix, and  $x(t):\RR_+ \to \RR^n$ be an arbitrary absolutely continuous function, also implying that $\dot x(t)$ exists almost everywhere.
Suppose that the following two conditions (particularizing \eqref{eq:general_condition} to $V(x) = \frac{1}{2} x^TAx$) holds for every $i=1,\dots,n$:
\begin{equation}\label{eq:good_sign}
 \dot x_i(t) A_{i,:}x(t)  \leq 0,
\end{equation}
where $ A_{i,:}$ denotes the $i$th row of $A$. Then
\begin{itemize}
\item[(a)] If $A$ is positive definite, $x(t)$ converges to a constant vector $x^*$.
\item[(b)] If $A$ is positive semi-definite and no nonzero vector of its kernel has a zero component ($w\in \ker A, w\neq 0 \Rightarrow w_i \neq 0, \forall i$), then either $x(t)$ converges to a constant vector $x^*$ or every accumulation point $\bar x$ of $x(t)$ lies in $\ker A$.
\end{itemize}
\end{theorem}
\medskip
Before presenting the proof, we note that condition (b) can only be satisfied if $\ker A$ has dimension 1. Indeed, if $v\neq w$ are linearly independent vectors in $\ker A$, one can always find a nontrivial linear combination $z = \alpha v + \beta w$ for which $z_i = 0$ for any given $i$. We also insist on $x$ being an arbitrary absolutely continuous function, with no assumption made about how it was generated. 

\medskip

\begin{IEEEproof}
We first show that (b) implies (a): Indeed, with $V(x):= \frac{1}{2} x^TAx$, it follows from \eqref{eq:good_sign} that $\dot V(x(t)) \leq 0$ so that $x(t)$ always remains in the set $\{x:x^T A x \leq V(0)\}$. When $A$ is positive definite, this set is compact and $x(t)$ has thus at least one accumulation point. Supposing that $x(t)$ would not converge, (b) implies that every accumulation point of $x(t)$ would be in the kernel of $A$, i.e. would be equal to 0, which implies that $x(t)$ would converge to 0 since it would be the only accumulation point. In $\RR^n$, convergence of a continuous trajectory is indeed equivalent to the existence of one single accumulation point. We therefore only need to prove (b) in the sequel.

Consider the hyperplane
\begin{equation}\label{eq:def_Ki}
K_i = \left\{x\in \RR^n~|~ A_{i,:} x = 0 \right\}
\end{equation}
orthogonal to the $i$th row $ A_{i,:}$ of $A$. If there is no accumulation point, or exactly one,
meaning that $x(t)$ converges to a constant vector, the statement (b) holds trivially. We suppose to the contrary that there exist multiple accumulation points. We select an arbitrary accumulation point $\bar x$ contained \emph{in the smallest possible number of hyperplanes $K_i$} and denote this smallest possible number by $k$. We will show that $\bar x\in \ker A$. Without loss of generality, we can assume the indices are ordered in such a way that
\begin{align*}
\bar{x} &\in K_1\cap K_2\cap \ldots \cap K_k,\\
\bar{x} &\notin K_{k+1}\cup K_{k+2}\cup \ldots \cup K_n.
\end{align*}
We can choose $\varepsilon >0$ such that two following two conditions hold: (i) $B(\bar{x},4\varepsilon) \cap
(K_{k+1}\cup K_{k+2}\cup \ldots \cup K_n) = \emptyset$  and (ii) there is at least one other accumulation point outside of $B(\bar{x},4\varepsilon)$ (otherwise $\bar x$ would be the only accumulation point). 
Due to the existence of this other accumulation point, $x(t)$ must infinitely often leave $B(\bar{x},3\varepsilon)$ while getting infinitely often into $B(\bar{x},\varepsilon)$. More precisely, there exists a diverging sequence of disjoint time intervals $[t^m_1,t^m_2]$ ($m=1,2,\dots$) such that $x(t^m_1) \in S(\bar{x},\varepsilon) = \partial B(\bar{x},\varepsilon)$,
$x([t^m_1,t^m_2])\subset cl~B(\bar{x},3\varepsilon)$ and $x(t^m_2) \in
S(\bar{x},3\varepsilon)= \partial B(\bar{x},3\varepsilon)$ for every $m$. We define $\Delta x^m = x(t^m_2) - x(t^m_1)$, so
$2\varepsilon \leq \norm{\Delta x^m} \leq 4 \varepsilon$. Our proof relies on the following two lemmas, which will establish that $(\Delta x^m)^TA\Delta x^m= \sum_i (\Delta x^m)_i(A\Delta x^m)_i\to 0$ as $m \to \infty$.

\begin{lemma}\label{lem:adx0i<k}
$\lim_{m\to \infty} (A \Delta x^m)_i = 0$ for $i\leq k$. As a consequence, $\lim_{m\to \infty} \Delta x^m_i (A \Delta x^m)_i = 0$ for $i\leq k$.
\end{lemma}
\begin{IEEEproof}
We first show that the distance between $x(t^m_2)$ and $K_i$ converges to 0 for every $i=1,\dots,k$. If it was not the case, there would be an infinite subsequence of $x(t^{m}_2)$ at a distance larger than some $\delta>0$ from $K_i$. Since the $x(t^m_2)$ are by definition in the compact set $B(\bar x, 3\varepsilon)$, this sequence would admit an accumulation point at a distance at least $\delta>0$ from  $K_i$. Moreover, this accumulation point could not belong to any $K_j$ with $j>k$ because these sets have no intersection with $B(\bar{x},4\varepsilon)$. Hence we would have an accumulation point that belongs to less than $k$ sets $K_i$, which contradicts the selection of $\bar x$ as an accumulation point of $x(t)$ belonging to the smallest possible number of $K_i$.

As a consequence the distance between $x(t^m_2)$ and every $K_i$, converges to 0, and a similar argument shows the same result for $x(t^m_1)$. This implies by definition of $K_i$ that $\lim_{m\to \infty}(A \Delta x^m)_i=\lim_{m\to \infty}A_{i:} x(t^m_2) - A_{i,:} x(t^m_1) = 0$ of $i\leq k$. The final implication of the Lemma follows from the boundedness of $\Delta x^m$.
\end{IEEEproof}

\smallskip

\begin{lemma}\label{lem:dx0i>k}
$\lim_{m\to \infty} (\Delta x^m)_i = 0$ for $i>k$. As a consequence, $\lim_{m\to \infty} \Delta x^m_i (A \Delta x^m)_i = 0$ for $i > k$.
\end{lemma}
\begin{IEEEproof}
For $i>k$, we know $B(\bar{x},3\varepsilon)$ is distant from
$K_i$ by at least $\varepsilon$. Hence, if $y\in
B(\bar{x},3\varepsilon)$, then $|A_{i,:}y|\ge c$ for some $c>0$. So since $x_i(t^m_2),x_i(t^m_1)\in B(\bar{x},3\varepsilon)$, we have
\begin{align*}
|x_i(t^m_2)-x_i(t^m_1)| &
   \le \int_{t^m_1}^{t^m_2} |\dot x_i(t)| dt \\
& \leq  \frac{1}{c}
\int_{t^m_1}^{t^m_2} |\dot x_i(t)| |A_{i,:}x(t)| dt.
\end{align*}
Condition \eqref{eq:good_sign} implies that $\dot x_i(t)$ and $A_{i,:}x(t)$ have opposite signs whenever they are both nonzero (and also for all other indices $j$),
hence we obtain from the previous inequality:
\begin{align*}
|x_i(t^m_2)-x_i(t^m_1)| &\leq - \frac{1}{c}\int_{t^m_1}^{t^m_2} \dot x_i(t)A_{i,:}x(t) dt \\
&\leq - \frac{1}{c}\int_{t^m_1}^{t^m_2} \sum_{j=1}^n\prt{\dot x_j(t)A_{j,:}x(t)} dt \\
& = \frac{1}{c}(V(x(t^m_1)) - V(x(t^m_2))),
\end{align*}
where we remind that $V(x) = \frac{1}{2}x^TAx$. This last inequality holds for every $m$, so that
$$
\sum_m |x_i(t_2^m)-x_i(t_1^m)| \le \frac{1}{c} \sum_m(V(x(t_1^m)) - V(x(t_2^m))) <
\infty,
$$
as $V(x(t))$ is non-increasing and the overall decrease of $V(x(t))$ is finite. 
Hence,
$|\Delta x_i^m| = |x_i(t_2^m)-x_i(t_1^m)| \rightarrow 0$ as $m\rightarrow \infty$ which we wanted to show. The last implication of the Lemma follows from the boundedness of $\Delta x^m$.
\end{IEEEproof}

\smallskip

It follows  from Lemmas \ref{lem:adx0i<k} and \ref{lem:dx0i>k} that
\begin{equation}\label{eq:DxADx->0}
\lim_{m\to \infty} (\Delta x^m)^T A \Delta x^m =0.
\end{equation}
We will now show that this implies that $k=n$, i.e. that $\bar x\in K_i$ for every $i$, and thus that $\bar x\in \ker A$ by the definition \eqref{eq:def_Ki} of the $K_i$. Suppose by contradiction that $k<n$, which implies that $\Delta x^m_n \to 0$ by Lemma \ref{lem:dx0i>k}. We claim that
\begin{equation}\label{eq:Dx-kernel >0}
  \liminf_{m\to\infty} dist(\Delta x^m, \ker A) \ge c > 0,
\end{equation}
for some 
$c$. Otherwise, since $4\varepsilon \ge \norm{\Delta x^m} \ge 2\varepsilon$,
an accumulation point of $\Delta x^m$ would reveal a vector $w \in \ker A$ with $4\varepsilon \ge ||w|| \ge 2\varepsilon$ and $w_n=0$ contradicting
our condition on the kernel in part (b) of the theorem statement. In turn, knowing that $A$ is positive semi-definite and \eqref{eq:Dx-kernel >0}, we get
$$
\liminf_{m\to\infty} (\Delta x^m)^T A (\Delta x^m) \ge c' > 0,
$$
for some positive $c'$. This is in contradiction with \eqref{eq:DxADx->0}.
Hence we must have $k=n$, meaning that $\bar x$ belongs to all $K_i$ and thus to $\ker A$.
Since $\bar x$ was selected as belonging to the smallest number of $K_i$ all others accumulation points also belong to all $K_i$ and thus to $\ker A$, which establishes the claim (b). This also implies claim (a) as explained in the first part of the proof.
\end{IEEEproof}

\smallskip

Observe that condition (b) of Theorem \ref{thm:elementwise} does not guarantee the existence of an accumulation point. And in case there is a single accumulation point, it may not be in $\ker A$, as the trajectory could for example stop anywhere (and thus converge)  without violating \eqref{eq:good_sign}. However, in case the trajectory has multiple accumulation points, they all belong to $\ker A$. It remains open to determine if (i) the condition on vectors with 0 entries in the kernel of $A$ can be relaxed, and (ii) if trajectories satisfying condition (b) may indeed diverge or have multiple accumulation points.

\subsection*{The particular case of Laplacian matrices}

We obtain stronger results for a specific class of positive semi-definite matrices: the (connected) graph Laplacians. A symmetric matrix $L$ is a Laplacian if all its off-diagonal entries are non-positive, i.e.  $L_{i,j}= - a_{ij} \le 0$ if $i\neq j$, and if each of its rows sums to 0, i.e. $L_{ii}  = \sum_{j\neq i} a_{ij}$ for all $i=1,\dots,n$. An $n\times n$ Laplacian is positive semi-definite, and has rank $n-1$ if the corresponding graph is connected, that is, every node can be reached from any other one in the graph defined by associating a node to each $i=1,\dots,n$ and connecting two nodes $i,j$ if $a_{ij}=a_{ji}>0$.
Laplacians play a major role in various disciplines, including algebraic graph theory \cite{biggs1993algebraic}, and are particularly important in consensus and synchronization applications, see e.g., \cite{franceschelli2013decentralized, wieland2011internal,patterson2014consensus}.

Laplacians have two properties of special interest in our context. First, observe that
\begin{equation}\label{eq:Lxi}
(Lx)_i = L_{ii}x_i  -\sum_{j\neq i} a_{ij} x_j = \sum _{j\neq i} a_{ij}(x_i-x_j),
\end{equation}
that is, $(Lx)_i$ is a weighted sum of the differences between $x_i$ and the other coordinates. Second,
$$
x^T L x = \sum_{i, j\neq i} a_{ij}(x_i-x_j)^2,
$$
i.e., the associated quadratic function is a weighted sum of the square differences between the $x_i$, and is thus a measure of the \quotes{disagreement} in $x$.
This also shows that the kernel of an $n \times n$ Laplacian is spanned by the vector $\1=[1,\dots,1]^T\in\R^n$, since the quadratic form above is 0 if and only if all $x_i$ are equal (recalling that the corresponding graph is connected).
We leverage these ideas to show that \eqref{eq:good_sign} implies convergence when the matrix is a Laplacian.
\medskip
\begin{theorem}\label{thm:laplacian}
Let $L$ be an $n\times n$ Laplacian whose corresponding graph is connected,
and let $x(t):\RR_+ \to \RR^n$ be an arbitrary absolutely continuous trajectory.
If
\begin{equation}\label{eq:good_sign_lap}
\dot x_i(t) (L x(t))_i \leq  0
\end{equation}
for every $i$, then $x(t)$ converges to a constant vector $x^*$,  and  $x_i(t)\in  [\min_j x_j(0), \max_j x_j(0)]$ for all $i,t$ so that
\begin{equation}\label{eq:convex_hull}
x^*_i \in [\min_j x_j(0), \max_j x_j(0)],~ \forall i.
\end{equation}
\end{theorem}
\medskip
\begin{IEEEproof}
We first show that $\min x_i(t)$ and $\max x_i(t)$ evolve monotonously. Note that the maximum of finitely many absolutely continuous functions is also absolutely continuous, and thus $\max_j x_j(t)$ has a derivative almost everywhere. Let $t$ be an arbitrary time at which this derivative and that of all the $x_i$ exists, and let  $I^*(t) = \{i:x_i(t)=\max_j x_j(t)\}$. Then we have from \eqref{eq:Lxi}
$$
(Lx(t))_i = \sum_{j} a_{ij}(x_i(t) - x_j(t)) \geq 0 \quad \forall i\in I^*(t),
$$
and \eqref{eq:good_sign_lap} implies for all $i\in I^*$ that $\dot x_i(t) \leq 0$.

By the continuity of all coordinates, there is a small enough $\varepsilon > 0$ such that for any $t-\varepsilon < t' < t + \varepsilon$ we have $\varnothing \neq I^*(t')\subseteq I^*(t)$, i.e., if $x_i(t) < \max_j x_j(t)$, then $x_i(t') < \max_j x_j(t')$ for all $t'$ in a small interval around $t$. This means that for any $|\delta| < \varepsilon$, there exists $i\in I^*(t) $ s.t.
$$
\frac{\max_j x_j(t+\delta)-\max_j x_j(t)}{\delta} = \frac{x_i(t+\delta)-x_i(t)}{\delta}.
$$
So when taking the limit $\delta \to 0$ we get $\frac{d}{dt} \max_j x_j(t) = \dot x_i(t)$ for one (or more) $i\in I^*(t)$, and we have seen that $\dot x_i(t) \leq 0$ for all $i\in I^*(t)$ so the same has to hold true for $\frac{d}{dt} \max_j x_j(t)$. Hence the absolutely continuous function $\max_j x_j(t)$ has a nonpositive derivative almost everywhere, which implies it is non-increasing. An analogous reasoning can be applied for $\min x_i(t)$.
As a consequence $x(t)$ always remains in the compact set $[\min_j x_j(0), \max_j x_j(0)]^n$ and has thus at least one accumulation point $\bar x$.

To argue by contradiction, assume now that
$x(t)$ does not converge. 
The kernel of $L$ is the set $\{\alpha \1\}$, and it follows thus from Theorem \ref{thm:elementwise} that the accumulation point $\bar x$ satisfies $\bar x = \bar \alpha \1$ for some $\bar \alpha$.
Since it
is an accumulation point, for every $\varepsilon$ there exist a time $t'$ at which $\varepsilon \geq | x_i(t') - \bar x_i|  = |x_i(t') - \bar \alpha|$ for every $i$. In particular, $\max_j x_j(t') \leq \bar \alpha + \varepsilon$ and $\min_j x_j(t') \geq \bar \alpha - \varepsilon$. The monotonicity of $\min_j x_j$ and $\max_j x_j$ implies then  $x_i(t) \in [\bar \alpha - \varepsilon, \bar \alpha + \varepsilon] $ for all $t>t'$. Since we can chose $\varepsilon$ arbitrarily small, $x(t)$ converges to $\bar x$,  contradicting our assumption. So $x(t)$  must indeed converge to some $x^*$, and the monotonicity of $\min_j x_j(t)$ and $\max_j x_j(t)$ implies \eqref{eq:convex_hull}.
\end{IEEEproof}

\section{Application to Platoons with bounded disturbances}\label{sec:platoon}

In this section, we study how to utilize condition \eqref{eq:general_condition} in designing a decentralized motion control scheme for the problem of keeping inter-agent distances in multi-vehicle-agent platoons at pre-defined desired values, using noisy inter-agent relative measurements, as considered in \cite{FidanAnderson07}. The paper \cite{FidanAnderson07} has proposed a deadzone based switching control scheme to solve this problem, guaranteeing to have the agent positions kept bounded, robustly to distance measurement noises with a known upper bound. In \cite{FidanAnderson07}, solution of the problem with the proposed control scheme is formally established only for two-agent platoons. Formal analysis for platoons with higher number of agents is left incomplete, ending with a conjecture on the agent positions being kept bounded and the inter-agent distances converging to certain intervals (balls) centered at the desired values, with radii proportional to the noise upper bound. The conjecture was supported by partial analysis for specific cases and simulation test results. The control scheme proposed in \cite{FidanAnderson07} is later adapted to the cooperative adaptive cruise control (CACC) problem of keeping a desired spacing between the consequent agents of a vehicle-platoon in \cite{SancarFidanHuissoon15}, introducing a moving frame of reference and considering the vehicle dynamics of the agents. Next, we revisit the problem considered in \cite{FidanAnderson07} in a more general setting to be defined in the following subsection, and propose an approach based on generation of agent trajectories satisfying the condition \eqref{eq:general_condition}.

\subsection{Problem}\label{sec:platoon_problem}

We consider a set of agents $1,\dots,n$ each with a position
$x_i(t) \in \RR$
A connected undirected graph $G$ represents the
possible sensing capabilities: ($(i,j)\in E$ implies that $i$ can
sense the relative position of $j$ with some noise, and
vice-versa). A particular case of graph is the \quotes{chain
graph}, with $E~=~\{(1,2),(2,3),\dots,(n-1,n)\}$.

The measures are subject to disturbance, so that if there is an edge $(i,j)\in E$ then agent $i$ can
sense $\hat \Delta_{ji}= x_j-x_i + w_{ji}$, where $w_{ji}$ is an arbitrary disturbance satisfying  $\abs{w_{ji}}\leq \bar w$, for some known $\bar w > 0$. The $w_{ji}$ are measurable, but not necessarily continuous.
For each
$(i,j)$ in $E$ we are given a desired distance $D_{ji}$, and the
ideal objective would be that for each $(i,j)$ , $x_j-x_i =
D_{ji}$. Those distances are supposed realizable, i.e., there exist
$p_1,\dots, p_n\in \RR$ such that $p_j-p_i =D_{ji}$ for all
$(i,j)\in E$. This implies in particular $D_{ij}=-D_{ji}$. This realizability constraint is automatically satisfied for the chain graph and for trees in general. For more general graphs, small mismatches of $D_{ij}$ could also be modeled as being part of the disturbances.

In the absence of communication between agents, it has been observed that
use of 
individual agent controllers in certain classical forms, such as proportional and proportional-integral, will lead to instabilities due to inconsistencies between the measurements of the inter-agent distance 
\cite{BaSu03,Bail06}.
Consider for example two agents 1, 2 with $D_{21}~=~-D_{12}~=~1$, and suppose $w_{21}=0.01$ while $w_{12}=0$, i.e.,
agent $1$ overestimates its distance to $2$. One can verify that if the agents use the same proportional controller based on the distance they sense, i.e., if each agent $i$ uses the control law $\dot x_i=\gamma (\hat{\Delta}_{ji}-D_{ji})$, where $j$
is the index of the other agent,  we will have $\dot x_1 + \dot x_2 = 0.01\gamma$, and hence the average position will move to infinity.
In the next subsection, we design a non-hierarchical control law for $\dot{x}_i(t)$  guaranteeing that all
$x_i$ remain bounded, and that all constraints are (asymptotically)
satisfied.

\subsection{Control Law}

For robustness to effects of the disturbances $w_{ij}$, we use non-linear threshold functions
such as
$$
T_w(x) = x \text{   if  } \abs{x}> w \text{   and   } 0 \text{  else}
$$
but any nondecreasing function for which $T_w(x)=0$ if and only if $\abs{x}\leq w$ can be used. These imply in particular that $xT_w(x+w')\geq 0$ for every $x$ if $\abs{w'}\leq w$.

The aim in our control law design is to have the agent move only when there is no doubt that it moves in the right direction.  For each agent, we propose the control law
\begin{equation}\label{eq:control_law_platoon}
\dot x_i = u_i = kT_{d(i)\bar w} \prt{\sum_{j|(i,j)\in E}
(\hat \Delta_{ji}-D_{ji})},
\end{equation}
where $d(i)$ is the degree of $i$ in the graph $G$.
Since $\Delta_{ij}$ differs from $x_i-x_j$ by at most $\bar w$, this control law implies that $u_i$ will be negative (resp. positive) if and only if $\sum_j (x_j-x_i) - D_{ji}$ is positive (resp. negative) for sure. We will show that \eqref{eq:control_law_platoon} guarantees convergence of $x$ to constant positions where the distance constraints are approximately satisfied, with errors that depend on $\bar w$ and the properties of the graph.

A similar control law was introduced independently in the context of consensus with unknown bounded disturbance in \cite{bauso2009consensus}. However, the final step of the convergence proof of  \cite{bauso2009consensus}, establishing convergence based on a condition akin to \eqref{eq:general_condition}
is inaccurate\footnote{Specifically, equation (18) in \cite{bauso2009consensus}, which the last arguments of the proof rely on, does not hold in general, indicating again the need for convergence results based on condition \eqref{eq:general_condition}.}.

\begin{figure}
\centering
\includegraphics[scale=0.43]{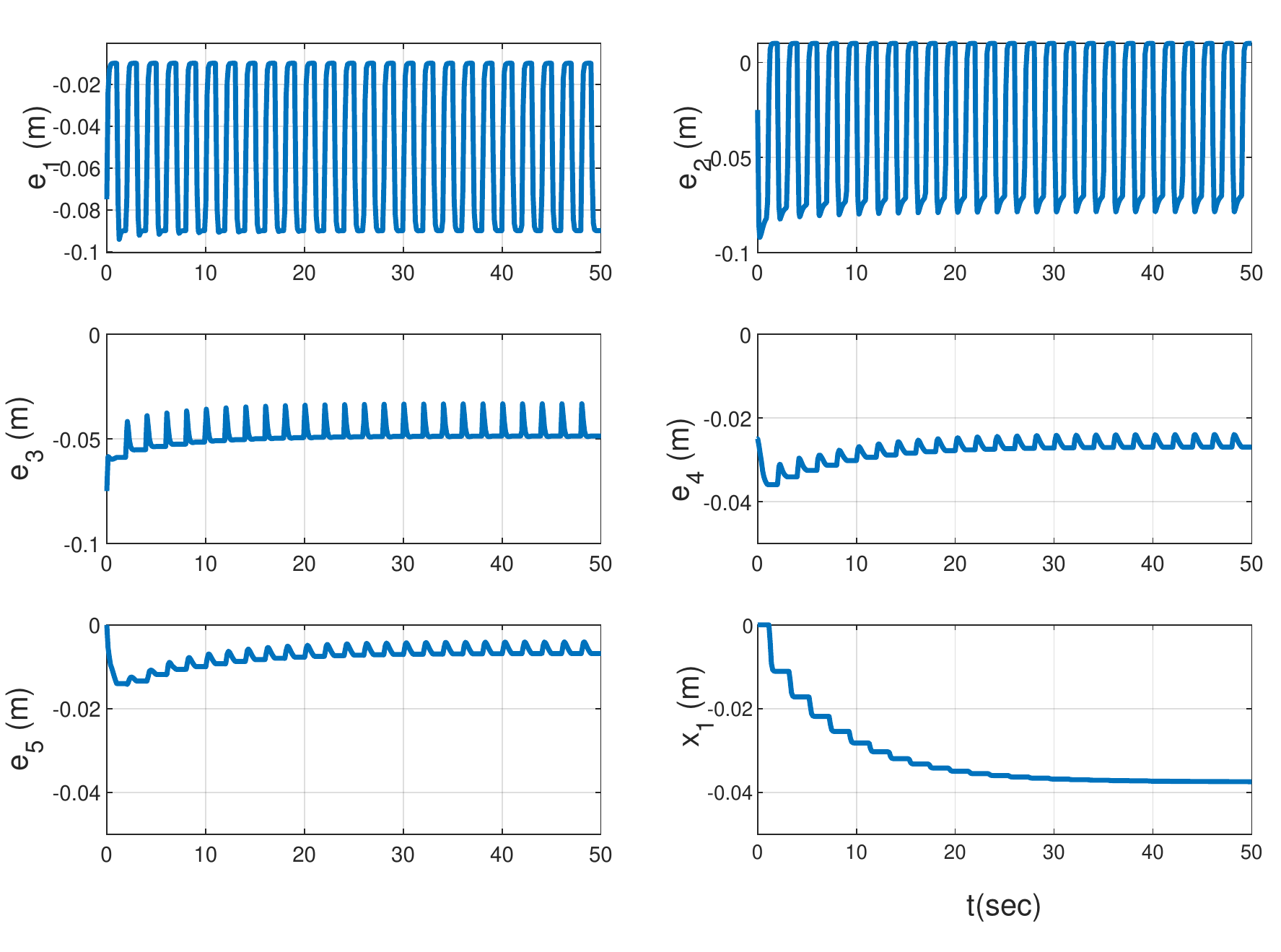}
\caption{\small{Example of evolution with time of the errors $e_i = x_{i+1}-x_i - D_{i,i+1}$ and of the position of the first agent $x_1$ for a 6 agent platoon with control laws \eqref{eq:deadzone_edge} and \eqref{eq:Tw2} and a chain sensing graph, showing that these control laws do not guarantee convergence.} }\label{fig:sim_MED07}
\end{figure}

\begin{figure}
\centering
\includegraphics[scale=0.43]{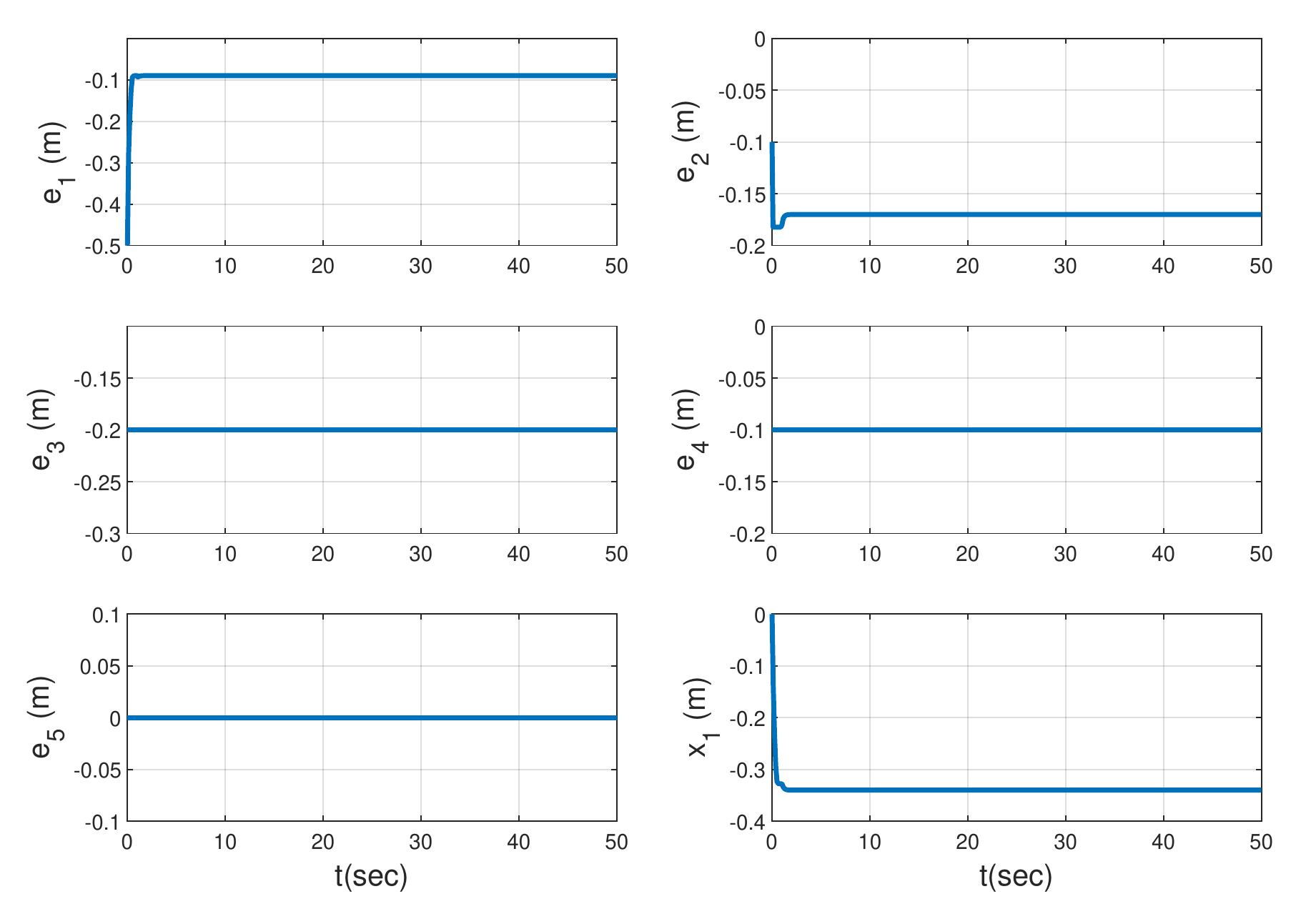}
\caption{\small{ Evolution with time of the errors $e_i = x_{i+1}-x_i - D_{i,i+1}$ and of the position of the first agent $x_1$ for a 6 agent platoon with control laws \eqref{eq:control_law_platoon} and \eqref{eq:Tw2} on the same sensing graph,  initial conditions and disturbances as in Fig. \ref{fig:sim_MED07}.}}\label{fig:sim_CDC19}
\end{figure}

We note that the issue of convergence is central here. For example, the similar looking control law
\begin{equation}\label{eq:deadzone_edge}
\dot x_i = u_i = k\sum_{j} T_{\bar w}(\hat \Delta_{ji}-D_{ji}),
\end{equation}
where the thresholds are applied to measurement as opposed to control actions, is observed to be inappropriate because agents would not necessarily converge to constant positions; for certain $w_{ij}$ they can indeed oscillate for ever. An example of such oscillations is presented in Fig. \ref{fig:sim_MED07} for a platoon with chain sensing graph with $n=6$ agents, where the initial positions  are $x_1(0)=0$, $x_2(0)=0.5$, $x_3(0)=1.4$, $x_4(0)=2.2$, $x_5(0)=3.1$, $x_6(0)=4.1$, the desired distances  $D_{21}=D_{32}=D_{43}=D_{54}=D_{65}=1$ (all in meter). For the sensor disturbances we take $w_{21}=w_{56}=0$, let
$w_{12}=w_{23}=w_{34}=w_{45}$ be a pulse signal with magnitude $0.1~m$, bias $-0.09~m$, period $2~$sec, and pulse width $1~$sec, and let
$w_{43}=w_{54}=w_{65}$ be a pulse signal with magnitude $0.1~m$,  bias $0.01~m$, period $2~$sec, and pulse width $1~$sec. Finally, $w_{32}$ is $1~$sec phase delayed version of $w_{43}$. The control law is that proposed in  \cite{FidanAnderson07}, i.e., \eqref{eq:deadzone_edge}, with
\begin{equation}
T_{\bar w}(x) =\left \{
\begin{array}{l}
x~\mbox{if}~\abs{x}> {\bar w} + \delta_w,\\
0~\mbox{if}~\abs{x}\leq \bar w,\\
(\abs{x}-\bar w) {\rm sgn}(x)/\delta_w~\mbox{else,}
\end{array}\right.
\label{eq:Tw2}
\end{equation}
$k=3$, $\bar w=0.1~$m, and $\delta_w=0.02~$m.
By comparison, Fig. \ref{fig:sim_CDC19} shows that the system converges when  control laws \eqref{eq:control_law_platoon} and \eqref{eq:Tw2} are used on the same initial conditions and disturbances.

\subsection{Convergence}

\begin{theorem}\label{thm:conv_platoon}
Consider $n$ agents $1,\dots,n$, with positions
$x_1(t),\dots,x_n(t) \in \RR$ at each time instant $t$, and a connected undirected sensing graph $G$ as detailed in Section \ref{sec:platoon_problem}. Under control law \eqref{eq:control_law_platoon}, for any realizable desired distances $D_{ij}$ (i.e. there exists $p_i$ such that  $D_{ji} = p_j-p_i,\forall (i,j)\in E$),
and any class of nondecreasing functions $T_w$ for which $T_w(x)=0$ if and only if $\abs{x}\leq w$,
\begin{itemize}
\item[(a)] $x(t)$ converges : $x^* = \lim _{t\to \infty}x(t)$ exists, and satisfies
\begin{equation}\label{eq:x*set_platoon}
\abs{\sum_{j:(i,j)\in E} (x_j^*-x_i^* -D_{ji})} \leq 2 d(i) \bar w,
\end{equation}
where $d(i)$ is the degree of agent $i$ in $G$ and $\bar w$ the bound on the disturbance.
\item[(b)] For every agent $i$ and all time $t$ there holds
\end{itemize}
\begin{equation}\label{eq:condition_inter_xp}
p_i + \min (x_j(0) - p_j)\leq x_i(t) \leq p_i + \max_j (x_j(0) - p_j),
\end{equation}
\end{theorem}
\begin{IEEEproof}
Let us perform a change of variable, defining $y_i = x_i - p_i$. Noting that $p_i-p_j =D_{ij}$, we have
\begin{align}
\dot y_i = \dot x_i &=  kT_{d(i)\bar w} \prt{\sum_{j|(i,j)\in E}
(y_j-y_i + w_{ji})}\nonumber
\\& = -kT_{d(i)\bar w} \prt{(Ly)_i - w_i}, \label{eq:control_y_compact}
\end{align}
where $w_i = \sum_{j|(i,j)\in E} w_{ji}$ satisfies, $\abs{w_i}\leq
d(i)\bar w$, and $L$ is the Laplacian matrix of the graph: $L_{ii}
= d(i)$, $L_{ji} = -1$ if $i$ and $j$ are connected, and 0 else. By definition of $T$ and in view of the bound $\abs{w_i}\leq
d(i)\bar w$, $\dot y_i$ can be positive only if  $(Ly)_i$ is negative, and vice versa, so that $\dot y( Ly)_i \leq 0$. Moreover, It is easy to confirm that $y$ is absolutely continuous as it is the integral of a measurable locally bounded function. Hence Theorem \ref{thm:laplacian} shows that $y$ converges to some $y^*$ and $y_i(t)$ remains at all time in $[\min y_i(t), \max y_i(t)]$, which implies the convergence of $x$ and the inclusion \eqref{eq:condition_inter_xp}.

We now prove that $\abs{(Ly^*)_i} \leq 2d(i) \bar w$, which implies \eqref{eq:x*set_platoon}, by contradiction.
Suppose this condition does not hold, and without loss of generality, that  $(Ly^*)_i > 2d(i) \bar w$. Since $y(t)$ converges to $y^*$, there is a time $t^*$ after which $(Ly(t))_i > 2d(i) \bar w + \alpha$ for some $\alpha >0$, and thus we have $(Ly(t))_i - w_i > d(i) \bar w + \alpha$, since  $w_i = \sum_{j|(i,j)\in E} w_{ij}$ satisfies, $\abs{w_i}\leq
d(i)\bar w$. Since $T$ is non-decreasing, this means there is a time after which $T_{d(i)\bar w}((Ly(t))_i - w_i)\geq T_{d(i)\bar w}(d(i)\bar w + \alpha)>0$, where the last inequality follows from $T_{d(i)\bar w}(z) = 0\Leftrightarrow \abs{z} \leq d(i) \bar w$.
As a result, $\dot y_i$ would remain negative and bounded away from 0 for all time $t>t^*$, in contradiction with its convergence to $y^*$. Hence we must have $\abs{(Ly^*)_i} \leq 2d(i) \bar w$ and thus \eqref{eq:x*set_platoon}.
\end{IEEEproof}
\smallskip

The convergence claim of Theorem \ref{thm:conv_platoon} remains valid if 
agents may stop for collision avoidance 
(see e.g. \cite{SancarFidanHuissoon15}) or for other reasons, as the inequality $\dot y_i (Ly)_i\leq 0$ used in the proof would still hold.

Theorem \ref{thm:conv_platoon} applies to any arbitrary connected sensing graph $G$. The particularization of \eqref{eq:x*set_platoon} to the line graph implies
$|x_2^*-x_1^*-D_{21}|  \leq 2 \bar w$
when applied to node 1, and
$
|2x_2^*-x_1^* - x_3^* - D_{21} - D_{32}|\leq 4\bar w.
$
when applied to node $2$, so that  $|x_3^*-x_2^* -D_{32}|  \leq 6\bar w$. An induction argument shows then
$$
|x_{\ell}^* - x_{\ell-1}^* -D_{\ell(\ell-1)}| \leq \min \prt{ 4\ell -6, 4n-4\ell-2} \bar w,
$$
where the second element in the min is obtained by starting the induction from the end of the platoon.
This illustrates, along the result (a) in Theorem  \ref{thm:conv_platoon}, that the control law  \eqref{eq:control_law_platoon} guarantees convergence of agents to constant positions
with the cost of having the upper bound of distance keeping errors dependent on the disturbance bound $\bar w$ and the number of agents.

In a similar setting, \cite{shi2019self} introduces a self-triggered scheme not relying on the knowledge of a bound on the disturbance, and guaranteeing bounded trajectories. For disturbances below a certain threshold they further show convergence by showing that $V=\frac{1}{2}x^TLx$ decreases by an amount that can be uniformly bounded from below after each triggering event, so that the total number of triggering times should be finite. 
Condition \eqref{eq:good_sign_lap} 
is satisfied in the framework of \cite{shi2019self}, so our results also directly ensure the convergence part of theirs.
This supports our hope that our results will serve 
as convenient \emph{tools} for stability analysis in different settings.

\section{Conclusion}\label{sec:ccl}

We have analyzed processes where the dynamics is not implicitly determined by (the gradient of) an energy function, but where that only serves as a barrier, leaving more freedom for the possible trajectory.
This framework beautifully matches the scenario of platoon formation, where the control of the dynamics has to be more conservative as it needs robustness as a priority over having an optimal configuration.

We have confirmed convergence of the processes when the energy function is quadratic described by a positive definite or Laplacian matrix. In a more general quadratic positive semi-definite case we have shown a partial concentration result, but we suspect much more is true.

These trajectory-based convergence results opens multiple perspectives: A straightforward challenge is to determine whether it is possible for a trajectory satisfying \eqref{eq:good_sign} to diverge or to have multiple accumulation points when $A$ has rank at most $n-1$ (and is not a Laplacian). Similarly, whether the absence of zero entries in the vectors of the kernel of $A$, required in condition (b) of Theorem \ref{thm:elementwise}, can be relaxed.
One obvious extension to broader context is to consider more general energy functions $V$ than quadratic ones.

Condition \eqref{eq:good_sign} can also be interpreted as requiring $\dot x$ and $\nabla V$ to belong to a same cone among a finite set of cone. This insight could be used to derive more general conditions with more general and/or position dependent cones.

\end{document}